\documentclass{amsart}

\usepackage{amsfonts,amscd,amsmath,latexsym}

\theoremstyle{plain}
\newtheorem{thm}{Theorem}[section]
\newtheorem{lemma}[thm]{Lemma}

\theoremstyle{remark}
\newtheorem{remark}[thm]{Remark}

\numberwithin{equation}{section}

\renewcommand{\a}{\alpha}
\newcommand{\p}{\mathfrak p}
\newcommand{\N}{\mathbb{N}}
\newcommand{\Z}{\mathbb{Z}}
\newcommand{\m}{\mathfrak m}
\newcommand{\Spec}{\operatorname{Spec}}
\newcommand{\sh}[1]{\ensuremath{{\mathcal #1}}}
\newcommand{\Proj}{\operatorname{Proj}}
\newcommand{\Hom}{\operatorname{Hom}}
\renewcommand{\mod}{\operatorname{mod}}
\newcommand{\Mod}{\operatorname{Mod}}
\newcommand{\End}{\operatorname{End}}
\newcommand{\proj}{\operatorname{proj}}
\newcommand{\GrSpec}{\operatorname{GrSpec}}
\newcommand{\QCoh}{\operatorname{QCoh}}
\newcommand{\Gr}{\operatorname{GrMod}}
\newcommand{\Tors}{\operatorname{Tors}}
\newcommand{\Coh}{\operatorname{Coh}}
\newcommand{\Inj}{\operatorname{Inj}}
\newcommand{\ass}{\operatorname{ass}}
\newcommand{\Kdim}{\operatorname{Kdim}}

\begin{document}
\title{The BGQ spectral sequence for noncommutative spaces}
\author{Christopher J. Pappacena}
\address{Department of Mathematics, Baylor University, Waco, TX 76798}
\email{Chris\_$\,$Pappacena@baylor.edu}
\thanks{The author was partially supported by a postdoctoral fellowhip from the Mathematical Sciences Research Institute, a summer sabbatical from Baylor University, and a grant from the National Security Agency.}

\begin{abstract} We prove an analogue of the Brown-Gersten-Quillen (BGQ) spectral sequence for noncommutative spaces.  As applications, we consider this spectral sequence over affine and projective spaces associated to right fully bounded noetherian (FBN) rings. \end{abstract} \maketitle

\section{Introduction}
One of the most fundamental results in the higher algebraic $K$-theory of schemes is the Brown-Gersten-Quillen (BGQ) spectral sequence \cite[Theorem 7.5.4]{Quillen}.  Recall that this spectral sequence is 
\begin{equation}
E_1^{pq}(X)=\bigoplus_{x\in X^p}K_{-p-q}(k(x))\Longrightarrow K_{n}'(X),\label{eq:BGQ}
\end{equation}
where $X$ is a noetherian scheme of finite Krull dimension.  Here $X^p$ is the set of codimension $p$ points of $X$, $k(x)$ is the residue field of $X$ at $x$, and $K_n'(X)$ is the $K$-theory of $\Coh(\sh{O}_X)$, the category of coherent $\sh{O}_X$-modules.   

As the title indicates, the goal of this note is to prove an analogue of this result for spaces studied in noncommutative algebraic geometry.  The main issues that arise are finding the appropriate noncommutative analogues of the geometric data that go into the spectral sequence \eqref{eq:BGQ}.  Our approach is to mimic the proof of the original formula (as found in \cite[section 7.5]{Quillen}) and let the spectral sequence tell us what the appropriate generalizations should be.  It will turn out that the appropriate noncommutative analogue of the underlying point set of $X$ will be the set of (isomorphism classes of) indecomposable injective $X$-modules $\Inj(X)$.  This gives further evidence that $\Inj(X)$ (called the \emph{injective spectrum} in \cite{Pappacena injective}) is an appropriate spectrum to study in noncommutative algebraic geometry.

After developing the noncommutative analogue of the BGQ spectral sequence in section 2, we apply our results to noncommutative spaces associated to right fully bounded noetherian (FBN) rings in section 3.  We show that if $R$ is right FBN, then the BGQ spectral sequence for $\Mod R$ can be described by data related to the prime ideals of $R$.  Similarly, if $R$ is a locally finite algebra over a field which is graded right FBN and generated in degree $1$, then we show that the BGQ spectral sequence for $\Proj R$ can be described in terms of the relevant graded prime ideals of $R$.

We briefly recall some of the definitions of noncommutative algebraic geometry.  A \emph{noncommutative space} $X$ is a Grothendieck category, also written as $\Mod X$ \cite{Smith subspaces}.  The objects of $\Mod X$ are called \emph{$X$-modules}.  The full subcategory of noetherian $X$-modules is written $\mod X$.  Examples of noncommutative spaces include $\QCoh(\sh{O}_X)$, the category of quasicoherent $\sh{O}_X$-modules when $X$ is a quasicompact and quasiseparated scheme,  $\Mod R$, the category of right modules over a ring $R$, and $\Proj R$, the noncommutative projective scheme associated to a finitely-generated graded $k$-algebra $R$ \cite{Artin Zhang}.  We refer the reader to \cite{Smith subspaces, Van blowup} for more details and information.

If $X$ is a noncommutative space, then as mentioned above the set of isomorphism classes of indecomposable injective $X$-modules is denoted $\Inj(X)$ and called the injective spectrum of $X$. For each $x\in\Inj(X)$, we fix a representative for its isomorphism class and denote it by $E(x)$.  If $E$ is an indecomposable injective $X$-module, then it is well-known that $\End(E)$ is a local ring, with maximal ideal (say) $J$.  We denote the division ring $\End(E)/J$ by $D(E)$; if $E\cong E(x)$ for $x\in\Inj(X)$ then we write $D(x)$ for $D(E)$ and call it the \emph{residue division ring at $x$}.

\subsection*{Acknowledgement} I would like to thank S. P. Smith for introducing me to noncommutative algebraic geometry, and for helping me formulate these results in their current degree of generality.

\section{The BGQ spectral sequence for a noncommutative space}
Let $X$ be a noetherian noncommutative space. By this we mean that $\Mod X$ is a locally noetherian category, so that every $X$-module is the direct limit of its noetherian submodules.  We fix an ordinal-valued dimension function $\dim$ on $\Mod X$ which is exact and finitely partitive, and set $\dim X=\sup\{\dim M:M\in\Mod X\}$ (provided the supremum exists). An example of such a dimension function is Krull dimension in the sense of Gabriel \cite[p. 382]{Gabriel}.  

Given an ordinal $\a$, and $X$-module $M$ is called \emph{$\a$-critical} if $\dim M=\a$ and $\dim M/N<\a$ for all nonzero submodules $N$ of $M$.  For $x\in\Inj(X)$, the \emph{critical dimension} of $x$ is the dimension of a critical submodule of $E(x)$ (which is well-defined but need not equal $\dim E(x)$ in general).

\begin{lemma} For each ordinal $\a$, let $\sh{T}_\a$ be the full subcategory of $\mod X$ consisting of those $X$-modules $M$ with $\dim M<\a$.

\begin{enumerate} 
\item[(a)] Each $\sh{T}_\a$ is a Serre subcategory of $\mod X$.
\item[(b)] If $S$ is a simple $\mod X/\sh{T}_\a$-module, then there exists an $\a$-critical $X$-module $M$ with $S\cong \pi_\a M$, where $\pi_\a:\mod X\rightarrow \mod X/\sh{T}_\a$ denotes the localization functor.
\end{enumerate}\label{simple lemma}
\end{lemma}

\begin{proof} Part (a) follows immediately from the fact that $\dim$ is exact, and (b) is a well-known consequence of the fact that $\dim$ is also finitely partitive (see for instance \cite[Proposition 3.9]{Pappacena injective}).
\end{proof}

From now on, assume that $X$ is a noetherian noncommutative space with $\dim X=d$ for some nonnegative integer $d$.  Then we have the following filtration of $\mod X$ by Serre subcategories:
\begin{equation}
0=\sh{T}_0\subset\sh{T}_1\subset\dots\subset \sh{T}_d\subset \sh{T}_{d+1}=\mod X.
\end{equation}
The hypothesis that $\dim$ is finitely partitive ensures that each object of the quotient category $\sh{T}_{i+1}/\sh{T}_i$ has finite length for all $i$.  Thus by \cite[Theorem 4, Corollary 1]{Quillen} we have isomorphisms for all $i$ and $n$:
\begin{equation}
K_n(\sh{T}_{i+1}/\sh{T}_i)\cong\bigoplus_{[S]\in\Sigma_i}K_n(\End(S)),\label{eq:Devissage}\end{equation}
where $\Sigma_i$ is the set of isomorphism classes of simple $\sh{T}_{i+1}/\sh{T}_i$-modules. For each $[S]\in\Sigma_i$, choose an $i$-critical $X$-module $M$ with $[S]=[\pi_iM]$. Then 
by \cite[Proposition 3.8]{Pappacena injective} we have that $\Sigma_i$ is in bijective correspondence with $\{x\in\Inj(X):\mbox{the critical dimension of $x$ is $i$}\}$.  Denote this latter set by $\Inj(X)_i$. Given $[\pi_iM]\in\Sigma_i$, there is a unique $x\in\Inj(X)_i$ with $E(M)\cong E(x)$.  By \cite[Proposition 4.1(b)]{Pappacena injective}, there is a ring isomorphism $\End(\pi_iM)\cong D(x)$.  Thus we can rewrite formula \eqref{eq:Devissage} as
\begin{equation}
K_n(\sh{T}_{i+1}/\sh{T}_i)\cong \bigoplus_{x\in\Inj(X)_i}K_n(D(x)).\label{eq:Devissage 2}
\end{equation}

Applying the Localization Theorem of $K$-theory \cite[Theorem 5]{Quillen} to each of the inclusions $\sh{T}_i\subset \sh{T}_{i+1}$ and using formula \eqref{eq:Devissage 2} gives long exact sequences
\begin{equation}
\cdots\longrightarrow K_n(\sh{T}_i)\longrightarrow K_n(\sh{T}_{i+1})\longrightarrow \bigoplus_{x\in\Inj(X)_i}K_n(D(x))\longrightarrow K_{n-1}(\sh{T}_i)\longrightarrow\cdots
\end{equation}

Splicing these long exact sequences together gives the desired generalization of the BGQ spectral sequence.
\begin{thm} Let $X$ be a noetherian noncommutative space, finite-dimensional with respect to a dimension function $\dim$.  Then there is a convergent spectral sequence
\begin{equation}
E_{pq}^1=\bigoplus_{x\in\Inj(X)_p}K_{p+q}(D(x))\Longrightarrow K_n(\mod X).\label{eq:NC BGQ}\end{equation}\label{BGQ}
\end{thm}

\begin{remark} (a) The spectral sequence \eqref{eq:NC BGQ} is of homological type because we filtered $\mod X$ by dimension instead of codimension.  Had we worked instead with a codimension function on $\Mod X$, then (under the same hypotheses as in Theorem \ref{BGQ}) we would have obtained a cohomological spectral sequence
\[E^{pq}_1=\bigoplus_{x\in\Inj(X)^p}K_{-p-q}(D(x))\Longrightarrow K_n(\mod X),\]
where $\Inj(X)^p$ is the set of indecomposable injectives of critical codimension $p$.

(b) If $X$ is a finite-dimensional noetherian scheme and we filter $\Coh(\sh{O}_X)$ by dimension of support, then the spectral sequence \eqref{eq:NC BGQ} reduces to the homological version of the usual BGQ spectral sequence \eqref{eq:BGQ}.  This follows by the bijective correspondence between indecomposable injective quasicoherent $\sh{O}_X$-modules and the underlying point set of $X$ \cite[Theorem VI.2.1]{Gabriel}.
\end{remark}

\section{Applications to right FBN rings}
\subsection{Affine spaces} In this subsection we consider the BGQ spectral sequence for the affine space $\Mod R$ when $R$ is a right FBN ring.  Recall that $R$ is said to be \emph{right FBN} if, for every prime ideal $\p$ of $R$, every essential right
ideal of $R/\p$ contains a nonzero two-sided ideal.  We shall use Krull dimension (in the sense of Gabriel) for our dimension function throughout.  Our goal is to describe the BGQ spectral sequence \eqref{eq:NC BGQ} for $\Mod R$ in terms of the prime ideals of $R$.  We begin with a lemma.

\begin{lemma} Let $R$ be a right noetherian ring, and let $U$ be a critical right ideal of $R/\p$.  Denote the injective hulls of $U$ in $\Mod R$ and $\Mod R/\p$ by $E$ and $E_\p$, respectively.  Then $D(E)\cong D(E_\p)$.\label{same D lemma}
\end{lemma}

\begin{proof} Note that we can identify $E_\p$ as an $R$-submodule of $E$.  Let $\tilde U$ denote the largest critical submodule of $E$, and similarly let $\tilde U_\p$ denote the largest critical submodule of $E_\p$ (so that $\tilde U_\p$ is a submodule of $\tilde U$).  
Then we have $D(E)\cong \Hom_R(U,\tilde U)$ \cite[Proposition 4.1(b)]{Pappacena injective}, and $D(E_\p)\cong\Hom_{R/\p}(U,\tilde U_\p)=\Hom_R(U,\tilde U_\p)$.  Now if $f\in\Hom_R(U,\tilde U)$, then $f(U)$ is annihilated by $\p$.  Hence $f(U)$ is in fact a submodule of $\tilde U_\p$, showing that $\Hom_R(U,\tilde U)=\Hom_R(U,\tilde U_\p)$.
\end{proof}

\begin{thm}  Let $R$ be a right FBN ring of finite Krull dimension. Then there is a convergent spectral sequence
\begin{equation}
E^1_{pq}=\bigoplus_{\p\in\Spec(R)_p}K_{p+q}(Q(R/\p))\Longrightarrow K_n(\mod
R), \label{eq:FBN BGQ}\end{equation}
where $\Spec(R)_p$ is the set of prime ideals of $R$ with $\Kdim(R/\p)=p$ and $Q(R/\p)$ is the Goldie quotient ring of $R/\p$.\label{FBNBGQ} \end{thm}

\begin{proof} Recall that the Gabriel correspondence is bijective for right fully bounded noetherian rings \cite[Theorem 8.13]{Goodearl Warfield}. Hence the map $x\mapsto\ass(E(x))$ gives a bijection between $\Inj(\Mod R)$ and $\Spec(R)$.  Moreover, the critical dimension of $E$ is exactly the Krull dimension of $R/\ass(E)$ \cite[Theorem 8.5]{GR}.  Thus there are bijections between $\Inj(\Mod R)_p$ and $\Spec(R)_p$ for all $p$.

Now let $\p\in\Spec(R)_p$, and let $E$ be an indecomposable injective $R$-module with $\ass(E)=\p$.  By the
proof of \cite[Theorem 8.5]{GR}, there is a uniform right ideal $U$ of $R/\p$
with $E\cong E(U)$, and without loss of generality we may choose $U$ to be $p$-critical.  Let $t$ be the uniform
dimension of $R/\p$, so that $U^{(t)}$ is essential in $R/\p$.  It follows
that $E(U)^{(t)}\cong E(U^{(t)})=E(R/\p)$.

Let $E_\p(U)$ denote the injective hull of $U$ in $R/\p$. Since $U$ is
a uniform right ideal of $R/\p$, we have that $E_\p(U)=QU$, where $Q=Q(R/\p)$.
Also, it is easy to see that $QU$ is itself a $p$-critical $R/\p$-module.  Now, $D(E_\p(U))\cong\End_{R/\p}(QU)\cong \End_R(QU)\cong D$, where $Q(R/\p)\cong M_t(D)$, and $D(E(U))\cong D(E_\p(U))$ by Lemma \ref{same D lemma}.  The fact that $K_n(D)\cong K_n(Q(R/\p))$ for all $n$ by Morita equivalence completes the proof. \end{proof}

 \subsection{Projective spaces} Suppose now that $R=\oplus_{n\geq 0} R_n$ is a right noetherian $\N$-graded $k$-algebra, $k$ a field. We shall assume that $R$ is locally finite, i.e. that $\dim_kR<\infty$ for all $n$. Our goal is to describe the BGQ spectral sequence for $\Proj R$ in the case where $R$ is graded right FBN. Here $\Proj R$ denotes the noncommutative projective scheme in the sense of Artin and Zhang \cite{Artin Zhang}. 

We begin by fixing notation and terminology. $\Gr R$ denotes the category of graded right $R$-modules, and we work exclusively in this category.  Thus all homomorphisms have degree 0, and $E(M)$ denotes the injective hull of $M$ in $\Gr R$.  We write $M(n)$ for the module obtained by shifting the grading on $M$ by $n\in\Z$.  We write $\GrSpec R$ for the set of homogeneous prime ideals of $R$, and call a prime ideal $\p$ \emph{relevant} if $\p$ does not contain $\m=\oplus_{n>0}R_n$.

As is customary, graded versions of notions in ring theory are defined using graded modules/homogeneous ideals.  Thus $R$ is graded right FBN if every graded essential right ideal of $R/\p$ contains a homogenous two-sided ideal, for every graded prime ideal $\p$ of $R$.  We refer the reader to \cite{NV} for further details and for any undefined terms.

If $R$ is graded right FBN, then the Gabriel correspondence $\Inj(\Gr R)\rightarrow \GrSpec R$ is no longer bijective.  However, it is proven in \cite[Theorem C.I.3.2]{NV} that if $R$ is graded right FBN and $E_1$ and $E_2$ are indecomposable injectives with $\ass(E_1)=\ass(E_2)$, then there is an integer $n$ such that $E_2\cong E_1(n)$. 

\begin{lemma} Let $R$ be a prime, graded right FBN $k$-algebra of finite Krull dimension, generated in degree $1$, and let $U$ be a graded critical right ideal of $R$ with injective envelope $E$.  Then there exists a positive integer $n$ such that $E\cong E(-n)$, and $E(R)\cong \oplus_{i=0}^{n-1}E(-i)^{(t_i)}$ for positive integers $t_0,\dots, t_{n-1}$.  \label{E(R) lemma}
\end{lemma}

\begin{proof}  Since $R$ is generated in degree $1$, there exists $x\in R_1$ such that $xU\neq 0$.  Note that $xU$ is a homogeneous right ideal of $R$, and right multiplication by $x$ induces a nonzero homomorphism $U(-1)\rightarrow xU$.  Since $U(-1)$ is critical and every homogeneous ideal of $R$ has the same Krull dimension, this homomorphism must be injective, showing that $xU\cong U(-1)$.  Thus $U(-1)$ is isomorphic to a homogeneous right ideal of $R$, and it follows that $U(-i)$ is isomorphic to a homogeneous right ideal of $R$ for all $i\geq 0$.  

Taking injective hulls we see that $E(-i)$ is isomorphic to a summand of $E(R)$ for all $i\geq 0$; since the decomposition of $E(R)$ into indecomposable injectives is unique up to isomorphism, there exist $i>j$ with $E(-i)\cong E(-j)$, so that $E(-(i-j))\cong E$.  Hence there is a smallest positive integer $n$ with the property that $E(-n)\cong E$, and the integers $t_0,\dots,t_{n-1}$ count the multiplicities with which $E,\dots,E(-(n-1))$ occur as summands of $E(R)$.
\end{proof}

\begin{lemma} Keeping the above notation, we have $\Hom(E(i),E(j))=0$ whenever $|j-i|<n$.\end{lemma}

\begin{proof} Twisting by $-i$, it suffices to show that $\Hom(E,E(j-i))=0$ for $|i-j|<n$.  Since $U$ is critical we have that every finitely-generated submodule of $E$ is critical. (The proof is a modification of \cite[Theorem 6.7]{GR} to the graded situation.) Thus $E$ is the sum of its critical submodules and so $E$ is critical \cite[Lemma 3.6]{Pappacena injective}. Thus any nonzero homomorphism $E\rightarrow E(j-i)$ must be injective, hence an isomorphism.  But $E\not\cong E(j-i)$ because $|i-j|<n$.
\end{proof}

Recall that if $R$ is graded prime, then the graded Goldie quotient ring $Q^{\rm gr}(R)$ of $R$ is graded simple, graded artinian by \cite[Theorem 4]{GS}. In particular the degree $0$ component $Q^{\rm gr}_0(R)$ is semisimple artinian.  The following lemma describes $Q^{\rm gr}_0(R)$.

\begin{lemma}  Keeping the above notation, we have $Q^{\rm gr}_0(R)\cong \oplus_{i=0}^{n-1} M_{t_i}(D(E))$.\label{endomorphism lemma}
\end{lemma}

\begin{proof} We have $Q^{\rm gr}(R)\cong E(R)$, so that $Q^{\rm gr}(R)\cong\oplus_{i=0}^{n-1}E(-i)^{(t_i)}$. Taking degree $0$ endomorphisms then gives
\[Q^{\rm gr}_0(R)\cong\End(Q^{\rm gr}(R))\cong\End(\oplus_{i=0}^{n-1}E(-i)^{(t_i)})\cong \oplus_{i=0}^{n-1}M_{t_i}(\End(E(-i)),\]
where the latter isomorphism follows from the previous lemma. Now $\End(E(-i))\cong \End(E)$ for all $i$ since twisting by $-i$ is an autoequivalence of $\Gr R$.  Since $\End(E)$ is local and $Q^{\rm gr}_0(R)$ is semisimple artinian we must have $\End(E)\cong D(E)$.
\end{proof}

We can now give a description of the BGQ spectral sequence for $\Proj R$ for certain graded right FBN $k$-algebras $R$, in terms of the relevant homogenous primes of $R$.  Recall that $\Proj R$ is defined to be the quotient category $\Gr R/\Tors R$, where $\Tors R$ is the full subcategory of torsion $R$-modules \cite[p. 233]{Artin Zhang}.  We denote by $\pi$ and $\omega$ the quotient and section functors associated to this localization.  Since $R$ is right noetherian and locally finite, $\Tors R$ consists precisely of those $M\in\Gr R$ with $\Kdim M=0$.  

\begin{thm} Let $R$ be a graded right FBN $k$-algebra of finite Krull dimension, generated in degree $1$.  Then there is a convergent spectral sequence
\begin{equation}
E^1_{pq}=\bigoplus_{\p\in\GrSpec(R)_{p+1}}K_{p+q}(Q^{\rm gr}_0(R/\p))\Longrightarrow K_n(\proj R),\end{equation}
where $\GrSpec(R)_i$ is the set of homogeneous prime ideals $\p$ with $\Kdim(R/\p)=i$.
\end{thm}

\begin{proof} Note that Krull dimension on $\Gr R$ induces Krull dimension on $\Proj R$, and that $\Kdim \sh{M}=\Kdim \omega\sh{M}-1$ for all $\sh{M}\in\Proj R$. If $\sh{E}$ is an indecomposable injective in $\Proj R$, then $\omega\sh{E}$ is an indecomposable injective in $\Gr R$ by \cite[Proposition 7.1]{Artin Zhang}. Moreover, if the critical dimension of $\sh{E}$ is $i$, then the critical dimension of $\omega\sh{E}$ is $i+1$.  Thus $\omega$ induces a bijection between $\Inj(\Proj R)_i$ and $\Inj(\Gr R)_{i+1}$ for all $i$.

So, the BGQ spectral sequence for $\Proj R$ can be written as 
\[E^1_{pq}=\bigoplus_{[E]\in\Inj(\Gr R)_{p+1}} K_{p+q}(D(E))\Longrightarrow K_n(\proj R).\]
Fix an element of $\Inj(\Gr R)_{p+1}$ with assassinator $\p$, and let $U$ be a graded critical right ideal of $R/\p$.  If $E$ is the injective hull of $U$, then $\Psi^{-1}(\p)$ consists of shifts of $E$, where $\Psi:\Inj(\Gr R)\rightarrow \GrSpec(R)$ is the Gabriel correspondence. Note that $E(i)\cong E(j)$ if and only if $E_\p(i)\cong E_\p(j)$, where $E_\p$ denotes the injective hull of $U$ in $\Gr R/\p$.  By Lemma \ref{E(R) lemma}, there is a positive integer $n_\p$ such that $\Psi^{-1}(\p)=\{E,E(-1),\dots,E(-(n_\p-1))\}$.

If we group the elements of $\Inj(\Gr R)_{p+1}$ together by their assassinators, and we use the fact that $D(E)\cong D(E(i))$ for all $i$, then we can rewrite the $E^1_{pq}$ term as $\bigoplus_{\p\in\GrSpec(R)_{p+1}}K_{p+q}(D(E))^{(n_\p)}$.  Now, $D(E)\cong D(E_\p)$ by a straightforward adaptation of Lemma \ref{same D lemma}, and $K_n(D(E_\p))^{(n_\p)}\cong K_n(Q^{\rm gr}_0(R/\p))$ for all $n$ by Morita equivalence and Lemma \ref{endomorphism lemma}.  Putting all of this together shows that we can write the $E^1_{pq}$ terms of the BGQ spectral sequence as $\bigoplus_{\p\in\GrSpec(R)_{p+1}}K_{p+q}(Q^{\rm gr}_0(R/\p))$, proving the theorem.
\end{proof}

\bibliographystyle{amsalpha}

 \end{document}